\documentclass[12pt,letterpaper]{article}
\usepackage{ifthen,latexsym,amssymb,amsmath}

\newcommand{\C}[1]{{\protect\cal #1}}
\newcommand{\B}[1]{{\bf #1}}
\newcommand{\I}[1]{{\mathbb #1}}

\newcommand{\ceil}[1]{\lceil #1\rceil}
\newcommand{\e}{\varepsilon}
\newcommand{\floor}[1]{\lfloor #1\rfloor}

\renewcommand{\mid}{:}

\newcommand{\beq}[1]{\begin{equation}\label{eq:#1}}
\newcommand{\eeq}{\end{equation}}
\newcommand{\req}[1]{\textrm{(\ref{eq:#1})}}

\newtheorem{theorem}{Theorem}
\newcommand{\bth}[2][nothing]{\ifthenelse{\equal{#1}{nothing}}
 {\begin{theorem}} {\begin{theorem}[#1]}\label{th:#2}}

\newtheorem{lemma}[theorem]{Lemma}
\newcommand{\blm}[2][nothing]{\ifthenelse{\equal{#1}{nothing}}
 {\begin{lemma}} {\begin{lemma}[#1]}\label{lm:#2}}

\newcommand{\bpf}[1][Proof.]{\smallskip\noindent{\it #1} }
\newcommand{\qed}{\nolinebreak\mbox{\hspace{5 true pt}%
  \rule[-0.85 true pt]{3.9 true pt}{8.1 true pt}}}
\newcommand{\epf}{\qed \medskip}

\newcommand{\brm}{\smallskip\noindent{\bf Remark.} }

\begin{document}

\renewcommand{\baselinestretch}{1.45}\rm

\renewcommand{\B}[1]{\mbox{\boldmath $#1$}}
\newcommand{\K}[1]{\C K_{#1}^{(k)}}
\newcommand{\h}[1]{H_{#1}^{(k)}}
\newcommand{\ex}{\mathrm{ex}}
\newcommand{\T}[1]{T^{(k)}(n,#1)}

\title{Exact Computation of the Hypergraph Tur\'an Function for
Expanded Complete $2$-Graphs} 
   
\author{Oleg Pikhurko\\
Department of Mathematical Sciences\\
Carnegie Mellon University\\
Pittsburgh, PA 15213\\
Web: {\tt http://www.math.cmu.edu/\symbol{126}pikhurko}}

\date{August 9, 2005}

\maketitle

\begin{abstract}
 Let $l>k\ge 3$. Let the $k$-graph $\h l$ be obtained from the
complete $2$-graph $K_l^{(2)}$ by enlarging each edge with a new set
of $k-2$ vertices. Mubayi [``A hypergraph extension of {Tur\'an's}
theorem'', to appear in \textit{J.\ Combin.\ Th.\ (B)}] computed
asymptotically the Tur\'an function $\ex(n,\h l)$. Here we determine
the exact value of $\ex(n,\h l)$ for all sufficiently
large $n$, settling a conjecture of Mubayi.\end{abstract}

\renewcommand{\baselinestretch}{1.45}\rm

\section{Introduction}

For $k,l\ge 2$ let $\K l$ be the family of all $k$-graphs $F$ with at
most ${l\choose 2}$ edges such that for some $l$-set $L$ (called the
\emph{core}) every pair $x,y\in L$ is covered by an edge of
$F$. Let the $k$-graph $\h l\in\K l$ be obtained from the complete
$2$-graph $K_l^{(2)}$ by enlarging each edge with a new set of $k-2$
vertices.

These $k$-graphs were recently studied by Mubayi~\cite{mubayi:04} in
the context of the \emph{Tur\'an $\ex$-function} which is defined as
follows. Let $\C F$ be a family of $k$-graphs. We say that a $k$-graph
$G$ is \emph{$\C F$-free} if no $F\in \C F$ is a subgraph of
$G$. (When we talk about subgraphs, we do not require them to be
induced.) Now, the \emph{Tur\'an function} $\ex(n,\C F)$ is the
maximum size of an $\C F$-free $k$-graph $G$ on $n$ vertices. Also,
let
 $$
 \pi(\C F)=\lim_{n\to\infty} \frac{\ex(n,\C F)}{{n\choose k}}.
 $$
 (The limit is known to exist, see Katona, Nemetz, and
Simonovits~\cite{katona+nemetz+simonovits:64}.)

To obtain the $k$-graph $\T l$, $l\ge k$, partition
$[n]=\{1,\dots,n\}$ into $l$ almost equal parts (that is, of sizes
$\floor{\frac nl}$ and $\ceil{\frac nl}$) and take those edges which
intersect every part in at most one vertex. Let us, for notational
convenience, identify $k$-graphs with their edge sets and, for a
$k$-graph $F$, write $\ex(n,F)$ for $\ex(n,\{F\})$, etc.

Mubayi~\cite[Theorem~1]{mubayi:04} proved the following result.

\bth[Mubayi]{M04} Let $n\ge l\ge k\ge 3$. Then $\ex(n,\K{l+1})=|\T
{l}|$, and $\T{l}$ is the unique maximum $\K{l+1}$-free $k$-graph of
order $n$.\qed\end{theorem}

It follows from Theorem~\ref{th:M04} and the super-saturation
technique of Erd\H os and Simonovits~\cite{erdos+simonovits:83} that
$\pi(\h l)=\pi(\K l)$, see~\cite[Theorem~2]{mubayi:04}. This gave us
the first example of a non-degenerate $k$-graph with known Tur\'an's
density for every $k$. (Previously, Frankl~\cite{frankl:90} did this
for all even $k$.) Settling a conjecture posed in~\cite{mubayi:04}, we
prove that the Tur\'an functions of $\h{l+1}$ and $\K{l+1}$ coincide
for all large $n$.

\bth{main} For any $l\ge k\ge 3$ there is $n_0(l,k)$ such that for any
$n\ge n_0(l,k)$ we have  $\ex(n,\h{l+1})=|\T
{l}|$, and $\T{l}$ is the unique maximum $\h{l+1}$-free $k$-graph of
order~$n$.\qed\end{theorem}

\brm Theorem~\ref{th:main} is true for $k=2$ by the Tur\'an
theorem~\cite{turan:41}. If $k\ge 3$ and $2\le l< k$, then
Theorem~\ref{th:main} is false:  $\ex(n,\K{l+1})=0$ while
$\ex(n,\h{l+1})>0$.\medskip 

\brm We do not compute an explicit upper bound on $n_0(l,k)$ as this
would considerably lengthen the paper. (For one thing, we would have
to reproduce some proofs from~\cite{mubayi:04} in order to calculate
an explicit dependence between the constants there.)\medskip

\section{Stability of $\h{l}$}

Two $k$-graphs $F$ and $G$ of the same order are \emph{$m$-close} if
we can add or remove at most $m$ edges from the first graph and make
it isomorphic to the second; in other words, for some bijection
$\sigma:V(F)\to V(G)$ the symmetric difference between
$\sigma(F)=\{\sigma(D)\mid D\in F\}$ and $G$ has at most $m$ edges.

Mubayi~\cite[Theorem~5]{mubayi:04} proved that $\K{l}$ is
\emph{stable}, meaning for the purpose of this article that for any
$\e>0$ there are $\delta>0$ and $n_0$ such that any $\K{l}$-free
$k$-graph $G$ of order $n\ge n_0$ and size at least
$(\pi(\K{l})-\delta){n\choose k}$ is $\e {n\choose k}$-close to
$\T{l-1}$. Here we prove the same statement for the single forbidden
graph $\h{l}$, which we will need in the proof of
Theorem~\ref{th:main}.

\blm{Hstable} For any $l> k\ge 3$ the $k$-graph $\h{l}$ is stable,
that is, for any $\e>0$ there are $\delta=\delta(k,l,\e)>0$ and
$n_0=n_0(k,l,\e)$ such that any $\h{l}$-free $k$-graph $G$ of order
$n\ge n_0$ and size at least $(\pi(\h{l})-\delta){n\choose k}$ is $\e
{n\choose k}$-close to $\T{l-1}$.\end{lemma}
 \bpf
 Let $\e>0$ be given. Choose $\delta>0$ which establishes the
stability of $\K l$ with respect to $\frac{\e}2$. Assume that
$\delta\le \e$. Let $n$ be large and $G$ be an $\h{l}$-free $k$-graph
on $[n]$ of size at least $(\pi(\h{l})-\frac{\delta}2) {n\choose k}$.

Let us call a pair $\{x,y\}$ of vertices \emph{sparse} if it is
covered by at most
 $$
 m=\left(l+(k-2){l\choose 2}\right){n\choose k-3}
 $$
 edges of $G$. Let $G'$ be obtained from $G$ by removing all edges
containing sparse pairs, at most ${n\choose 2}\times m<
\frac{\delta}{2} {n\choose k}$ edges.

Let us show that the $k$-graph $G'$ is $\K l$-free. Suppose on the
contrary that every pair from some $l$-set $L$ is covered by an edge
of $G'$. It follows that every pair $\{x,y\}\subset L$ is not sparse
with respect to $G$, that is, $G$ has more than $m$ edges containing
$\{x,y\}$. This means that if we have a partial embedding of $\h l$
into $G$ with the core $L$, then we can always find a $G$-edge $D\ni x,y$
such that $D\setminus\{x,y\}$ is disjoint from the rest of the
embedding. Thus $G$ has an $\h l$-subgraph with the core $L$, a
contradiction.

We have $|G'|\ge (\pi(\h{l})-\delta) {n\choose k}$. By the stability
of $\K l$, $G'$ is $\frac{\e}2\, {n\choose k}$-close to $\T
{l-1}$. The triangle inequality implies that $G$ is
$(\frac{\delta}2+\frac{\e}2){n\choose k}$-close to $\T {l-1}$. As
$\delta\le \e$, this finishes the proof of the lemma.\epf

\section{Exactness}

\noindent\textit{Proof of Theorem~\ref{th:main}.}  Let us choose, in
this order, positive constants $c_1,\dots,c_5$, each being
sufficiently small depending on the previous constants. Then, let
$n_0$ be sufficiently large. In fact, we can take some simple explicit
functions of $k,l$ for $c_1,\dots,c_5$. However, $n_0$ should also be
at least as large as the function $n_0(k,l+1,c_5)$ given by
Lemma~\ref{lm:Hstable}.

Let $G$ be a maximum $\h {l+1}$-free graph on $[n]$ with $n\ge
n_0$. We have
 \beq{g>t}
 |G|\ge |\T l|\ge \frac{l\,(l-1)\dots
   (l-k+1)}{l^k}\, {n\choose k}=\pi(\h{l+1})\, {n\choose k},
 \eeq
 where the first inequality follows from the fact that $\T l$ is
$\h{l+1}$-free while the second inequality can be shown directly. (For
example, a simple averaging shows that the function $ |\T l|/{n\choose
k}$ is decreasing in $n$.)

Let $V_1\cup \dots\cup V_l$ be a partition of $[n]$ such that
 $$
 f
 =   \sum_{D\in G}\, \Big| \{ i\in[l]\mid D\cap
V_i\not=\emptyset\}\Big|
 $$
 is maximum possible. Let $T$ be the complete $l$-partite $k$-graph on
$V_1\cup \dots\cup V_l$. Clearly, $f\ge k\,|T\cap G|$. As $n$ is
sufficiently large, Lemma~\ref{lm:Hstable} implies that $G$
is $c_5{n\choose k}$-close to $\T{l}$. (The value of $\delta>0$
returned by Lemma~\ref{lm:Hstable} is not significant here because of
the lower bound~\req{g>t} on the size of $G$.)  The choice of $T$
implies that $f\ge k(|G|-c_5{n\choose k})$. On the other hand, $f\le
k|G|-|G\setminus T|$. It follows that
 \beq{GminusT}
 |G\setminus T|\le c_5 k {n\choose k}.
 \eeq
 
Thus we have $|T|\ge|\T l|- c_5 k {n\choose k}$.  This bound on $|T|$
can be easily shown to imply (or, alternatively, see Claim~1
in~\cite[Proof of Theorem~5]{mubayi:04}) that for each $i\in[l]$ we
have, for example,
 \beq{2}
 |V_i|\ge \frac n{2l}.
 \eeq

Let us call the edges in $T\setminus G$ \emph{missing} and the edges
in $G\setminus T$ \emph{bad}. As $|T|\le |\T l|$ with equality if and
only if $T$ is isomorphic to $\T l$,
see~\cite[Equation~(1)]{mubayi:04}, the number of bad edges is at
least the number of missing edges. It also follows that if $G\subset
T$, then we are done. Thus, let us assume that $B$ is non-empty, where
the $2$-graph $B$ consists of all \emph{bad} pairs, that is, pairs of
vertices which come from the same part $V_i$ and are covered by an
edge of $G$.

For vertices $x,y$ coming from two different parts $V_i$, call the
pair $\{x,y\}$ \emph{sparse} if $G$ has at most 
 $$
 m=\left({l+1\choose
2}(k-2)+l+1\right){n\choose k-3}
 $$
 edges containing both $x$ and $y$; otherwise $\{x,y\}$ is called
\emph{dense}.

Note that there are less than $c_4 n^2$ sparse pairs for otherwise we
get a contradiction to~\req{GminusT}: each sparse pair generates at
least 
 \beq{sparse}
 \left(\frac{n}{2l}\right)^{k-2}- m \ge \frac12\,
\left(\frac{n}{2l}\right)^{k-2}
 \eeq
 missing edges by~\req{2} while each missing edge contains at most
${k\choose 2}$ sparse pairs.

Take any bad pair $\{x_0,x_1\}$, where, for example, $x_0,x_1\in V_1$
are covered by $D\in G$.  The number of vertices in $\h{l+1}$ is
${l+1\choose 2}(k-2)+l+1$. Therefore, if we have a partial
embedding of $\h{l+1}$ into $G$ such that a pair of vertices $x,y$
from the core is dense, then we can find a $G$-edge containing both
$x,y$ and disjoint from the rest of the embedding. It follows that for
any choice of $(x_2,\dots,x_l)$, where $x_i\in V_i\setminus D$ for
$2\le i\le l$, at least one pair $\{x_i,x_j\}$ with
$\{i,j\}\not=\{0,1\}$ is sparse. Since $x_0$ and $x_1$ are fixed, each
such sparse pair $\{x_i,x_j\}$ is counted, very roughly, at most
$n^{l-3}$ times if $\{x_i,x_j\}\cap \{x_0,x_1\}=\emptyset$, and at
most $n^{l-2}$ times if $\{x_i,x_j\}\cap \{x_0,x_1\}\not=\emptyset$.

Since we have at most $c_4 n^2$ sparse pairs, the number of times the
former alternative occurs is at most 
 $$
 c_4 n^2 \times n^{l-3}\le \frac 12
\left(\frac{n}{2l}-k\right)^{l-1}.
 $$
 That is, by~\req{2}, for at least half of the choices of
$(x_2,\dots,x_l)$, the obtained sparse pair intersects
$\{x_0,x_1\}$. Let $A$ consist of those $z\in V(G)$ which are incident
to at least $c_1n$ sparse pairs. Since $\frac14 (\frac{n}{2l}-k)^{l-1}
/ n^{l-2} \ge c_1 n$, at least one of $x_0$ and $x_1$
belongs to $A$. Thus, in summary, we have proved that every bad pair
intersects $A$.

Considering the sparse pairs, we obtain by~\req{sparse} at least 
 $$
 \frac{|A|\times c_1n}2 \,\times \frac12\,\left(\frac{n}{2l}\right)^{k-2} \times
{k\choose 2}^{-1}\ge |A|\times c_2n^{k-1}.
 $$
 missing edges and, consequently, at least $|A|\times c_2 n^{k-1}$ bad
edges. Let $\C B$ consist of the pairs $(D,\{x,y\})$, where
$\{x,y\}\in B$, $D\in G$ and $x,y\in D$. (Thus $D$ is a bad edge.) As
each bad edge contains at least one bad pair, we conclude that $|\C
B|\ge |A|\times c_2 n^{k-1}$. For any $(D,\{x,y\})\in \C B$, we have
$\{x,y\}\cap A\not=\emptyset$. If we fix $x$ and $D$, then, obviously,
there are at most $k-1$ ways to choose a bad pair $\{x,y\}\subset
D$. Hence, some vertex $x\in A$, say $x\in V_1$, belongs to at least
 \beq{c2}
 \frac{|\C B|}{(k-1)\,|A|} \ge \frac{c_2}{k-1}\, n^{k-1}
 \eeq
 bad edges, each intersecting $V_1$ in another vertex $y$.

Let $Y\subset V_1$ be the neighborhood of $x$ in the $2$-graph $B$. We
have 
 $$
 |Y|\ge \frac{c_2}{k-1}\, n^{k-1} \times {n\choose k-2}^{-1}\ge c_3n.
 $$
 For $j\in[2,l]$ let $Z_j$ consist of those $z\in V_j$ for which
$\{x,z\}$ is dense.

Suppose first that $|Z_j|\ge c_3n$ for each $j\in[2,l]$. In this case
we do the following. For every $y\in Y$, fix some $D_{y}\in G$
containing both $x$ and $y$. Consider an $(l+1)$-tuple
$L=(x,y,z_2,z_3,\dots,z_l)$, where $y\in Y$ and $z_j\in Z_j\setminus
D_{y}$ are arbitrary. We can find a partial embedding of $\h{l+1}$
with core $L$ such that every pair containing $x$ is covered: the pair
$\{x,y\}$ is covered by $D_{y}$ while each pair $\{x,z_i\}$ is
dense. Since $G$ is $\h{l+1}$-free, at least one pair from the set
$\{y,z_2,\dots,z_l\}$ is sparse. Since there are at least $(c_3
n-k)^l$ choices of $L$ (note that $x$ is fixed), this gives us at least
$(c_3 n-k)^l/n^{l-2}> c_4n^2$ sparse pairs, which is a contradiction
as we already know.

Hence, assume that, for example, $|Z_2|< c_3n$. This means that all
but at most $c_3 n$ pairs $\{x,z\}$ with $z\in V_2$ are sparse, that
is, there are at most 
 \beq{xV2}
 c_3n\times {n\choose k-2}+ n\times m\le c_3n^{k-1}
 \eeq
 $G$-edges containing $x$ and intersecting $V_2$. Let us contemplate
moving $x$ from $V_1$ to $V_2$. Some edges of $G$ may decrease their
contribution to $f$ by $1$. But each such edge must contain $x$ and
intersect $V_2$ so the corresponding total decrease is at most $c_3
n^{k-1}$ by~\req{xV2}. On the other hand, the number of edges of $G$
containing $x$, intersecting $V_1\setminus\{x\}$, and disjoint from
$V_2$ is at least $\frac{c_2}{k-1}\,n^{k-1}-c_3 n^{k-1}$ by~\req{c2}
and~\req{xV2}. As $c_3$ is much smaller than $c_2$, we strictly
increase $f$ by moving $x$ from $V_1$ to $V_2$, a contradiction to the
choice of the parts $V_i$. The theorem is proved.\epf

\section{Concluding Remarks}

Lemma~\ref{lm:Hstable} also follows from the following more general
Lemma~\ref{lm:stability}. In order to state the latter result, we need
some further definitions.

Let us call a family $\C F$ of $k$-graphs \emph{$s$-stable} if for any
$\e>0$ there are $\delta>0$ and $n_0$ such that for arbitrary $\C
F$-free $k$-graphs $G_1,\dots,G_{s+1}$ of the same order $n\ge n_0$,
each of size at least $(\pi(\C F)-\delta){n\choose k}$, some two are
$\e {n\choose k}$-close. Please note that if $\C F$ is $s$-stable for
some $s$ then it is also $t$-stable for any
$t>s$. Lemma~\ref{lm:Hstable} implies that $\h{l}$ is $1$-stable.  Let
$F[t]$ denote the \emph{$t$-blowup} of a $k$-graph $F$, where each
vertex $x$ is replaced by $t$ new vertices and each edge is replaced
by the corresponding complete $k$-partite $k$-graph.  Clearly,
$|F[t]|=t^k\,|F|$.

\blm{stability} Let $t\in\I N$. Let $\C F$ be a finite family of
$k$-graphs which is $s$-stable. Let $\C H$ be another (possibly
infinite) $k$-graph family such that for each $F\in \C F$ there is
$H\in \C H$ such that $H\subset F[t]$. If $\pi(\C H)\ge \pi(\C F)$,
then $\pi(\C H)= \pi(\C F)$ and $\C H$ is $s$-stable.\end{lemma}
 \bpf Our proof uses the following theorem of R\"odl and
Skokan~\cite[Theorem~7.1]{rodl+skokan:04u} which in turn relies on the
Hypergraph Regularity Lemma of R\"odl and Skokan~\cite{rodl+skokan:04}
and the Counting Lemma of Nagle, R\"odl, and
Schacht~\cite{nagle+rodl+schacht:04} (see also Gowers~\cite{gowers:05}).

\bth[R\"odl and Skokan]{rsF} For all integers $l>k\ge 2$ and a real $\e>0$
there exist $\mu=\mu(k,l,\e)>0$ and $n_1=n_1(k,l,\e)\in\I N$
such that the following statement holds.

Given a $k$-graph $F$ with $v\le l$ vertices, suppose that a $k$-graph
$G$ with $n>n_1$ vertices contains at most $\mu n^v$ copies of $F$ as
a subgraph. Then one can delete at most $\e {n\choose k}$ edges of $G$
to make it $F$-free.\qed\end{theorem}

Let $\e>0$ be arbitrary. Let $\delta>0$ and $n_0$ be constants
satisfying the $s$-stability assumptions for $\C F$ and
$\frac{\e}3$. Assume that $\delta\le \e$. Let $l$ be the maximum order
of a $k$-graph in $\C F$ and $m=|\C F|$.  Let
$\mu=\mu(k,l,\frac{\delta}{3m})$ and $n_1=n_1(k,l,\frac{\delta}{3m})$
be given by Theorem~\ref{th:rsF}. Also, assume that $n_2$ is so large
that for every $F\in \C F$ any $F[t]$-free $k$-graph of order $n\ge
n_2$ contains at most $\mu n^{v(F)}$ copies of $F$, where $v(F)$
denotes the number of vertices in $F$. Such $n_2$ exists because any
$F[t]$-free $k$-graph $G$ of order $n$ has at most $o(n^{v(F)})$
copies of $F$, which follows from a theorem of Erd\H os~\cite{erdos:64}. Let
$n_3=\max(n_0,n_1,n_2)$.

Let $n\ge n_3$ and let $G_1,\dots,G_{s+1}$ be arbitrary $\C H$-free
$k$-graphs each having $n$ vertices and at least $(\pi(\C
F)-\frac{\delta}2){n\choose k}$ edges. By Theorem~\ref{th:rsF} (and the
choice of $n_1$ and $n_2$), for each $F\in
\C F$ each $G_i$ can be made $F$-free by removing at most
$\frac{\delta}{3m} {n\choose k}$ edges. Hence, we can transform $G_i$
into an $\C F$-free $k$-graph $G_i'\subset G_i$ by removing at most
$|\C F|\,\frac{\delta}{3m} {n\choose k} \le \frac{\delta}{3} {n\choose
k}$ edges.

We conclude that $\pi(\C F)\ge \pi(\C H)- \frac{\e}{3}$. As $\e>0$ was
arbitrary, we have $\pi(\C F)=\pi(\C H)$. Thus the density of each
$G_i'$ is at least $\pi(\C H)-\frac{\delta}2- \frac{\delta}3 > \pi(\C F)-
\delta$. By the $s$-stability of $\C F$, some two of these graphs, for
example, $G_i'$ and $G_j'$, are $\frac{\e}3 {n\choose 2}$-close. It
follows that $G_i$ and $G_j$ are $\e {n\choose k}$-close. Thus the
constants $\frac{\delta}2$ and $n_3$ demonstrate the $s$-stability of
$\C H$, proving Lemma~\ref{lm:stability}.\epf

The line of argument we used in this article might be useful for
computing the exact value of $\ex(n,F)$ for other forbidden $k$-graphs
$F$. The approach in general could be the following.
 \begin{enumerate}
 \item Find a suitable $k$-graph family $\C F\ni F$ for which we can
compute $\pi(\C F)$ and prove the stability of $\C F$.
 \item Deduce from Lemma~\ref{lm:stability} that $\pi(F)=\pi(\C F)$
and $F$ is stable too.
 \item Using the stability, obtain the exact value of $\ex(n,F)$. (The
fact that stability often helps in proving exact results for the
hypergraph Tur\'an problem was observed and used by F\"uredi and
Simonovits~\cite{furedi+simonovits:03}, Keevash and
Sudakov~\cite{keevash+sudakov:03b,keevash+sudakov:03a}, and others.)
 \end{enumerate}

Extending the results by Sidorenko~\cite{sidorenko:87}, the
author~\cite{pikhurko:05c} has successfully applied the above
approach to computing the exact value of $\ex(n,T^{(4)})$ for $n\ge
n_0$, where the $k$-graph $T^{(k)}$ consists of the following three
edges: $[k]$, $[2,k+1]$, and $\{1\}\cup [k+1,2k-1]$. The exact value
of $\ex(n,T^{(3)})$ was previously computed by Frankl and
F\"uredi~\cite{frankl+furedi:83} (see also
Bollob\'as~\cite{bollobas:74}, Keevash and
Mubayi~\cite{keevash+mubayi:04}).

Lemma~\ref{lm:Hstable} has an interesting application.  Namely, the
method of Mubayi and the author~\cite{mubayi+pikhurko:03} (combined
with Lemma~\ref{lm:Hstable}) shows that the pair
$(\h{k+2},K_{k+1}^{(k)})$ is \emph{non-principal} for any $k\ge 3$,
that is,
 \beq{np}
 \pi\left(\{\h{k+2},K_{k+1}^{(k)}\}\right) <
\min\left\{\pi(\h{k+2}),\pi(K_{k+1}^{(k)})\right\},
 \eeq
 where $K_m^{(k)}$ denotes the complete $k$-graph of order $m$. This
completely answers a question of Mubayi and
R\"odl~\cite{mubayi+rodl:02} (cf.\ also Balogh~\cite{balogh:02}). We
refer the Reader to~\cite{mubayi+pikhurko:03} for further details.

\section*{Acknowledgments}

The author is grateful to Vojta R\"odl and Mathias Schacht for
providing the manuscripts~\cite{nagle+rodl+schacht:04,rodl+skokan:04u}
before their publication and to the anonymous referees for the very
useful and detailed comments.

\renewcommand{\baselinestretch}{1.0}

\small

\providecommand{\bysame}{\leavevmode\hbox to3em{\hrulefill}\thinspace}
\providecommand{\MR}{\relax\ifhmode\unskip\space\fi MR }
\providecommand{\MRhref}[2]{%
  \href{http://www.ams.org/mathscinet-getitem?mr=#1}{#2}
}
\providecommand{\href}[2]{#2}

\end{document}